\documentstyle[12pt]{article}

\newcommand{\be}{\begin{equation}}
\newcommand{\ee}{\end{equation}}
\newcommand{\zz} {\overline{z}}
\newcommand{\nn} {\noindent}
\newtheorem{proposition}{Proposition}[section]
\newtheorem{definition}{Definition}[section]
\begin{document}
\title{Deformed Harmonic Oscillator Algebras defined by their Bargmann representations.}
\author{Mich\`{e}le IRAC-ASTAUD and Guy RIDEAU\\ 
Laboratoire de Physique Th\'{e}orique de la mati\`{e}re condens\'ee\\
Universit\'{e} Paris VII\\2 place Jussieu F-75251 Paris Cedex 05, FRANCE}
\date{}
\maketitle
 to appear in Reviews in Mathematical Physics.

\begin{abstract}
Deformed Harmonic Oscillator Algebras are generated by  four operators 
two mutually adjoint 
 $a$ and $ a^\dagger$, and two self-adjoint $N$ and the unity $1$
  such as : 
  
  \nn $[a,N] = a ,\quad [a^\dagger ,N]=
-a^\dagger $, $a^\dagger a = \psi(N)$ and $\quad aa^\dagger =\psi(N+1)$.

\nn The Bargmann Hilbert space  is defined as a space of functions, holomorphic in
  a ring of the complex plane, equipped with a scalar product involving
   a true integral. In a Bargmann representation, the operators of a 
    Deformed Harmonic Oscillator Algebra act on a  Bargmann Hilbert space and
	the creation (or the annihilation operator) is the multiplication by $z$.
 We discuss the conditions of existence of Deformed Harmonic Oscillator Algebras
  assumed
 to admit a given Bargmann representation.
\end{abstract}

\section{Introduction}

\nn In previous papers  \cite{canada}, \cite{barg1}, \cite{barg2},
 we introduced what we have called
  Deformed Harmonic Oscillator Algebras (DHOA).

\begin{definition} \label{def:DHOA}
A Deformed Harmonic Oscillator Algebra is a free algebra generated by  four 
operators :
 
\nn - the annihilation operator $a$,
 the creation operator $a^\dagger$ that are mutually adjoint, 

\nn - the  self-adjoint energy operator
  $N$ and the unity $1$ 

\nn  satisfying the following commutation relations :

\begin{equation}
[a,N] = a ,\quad [a^\dagger ,N] =
-a^\dagger ,\quad
a^\dagger a = \psi(N) , \quad aa^\dagger =\psi(N+1)
\label{a3}
\ee

\nn where $\psi$ is a real analytical function. 

\end{definition}

\nn When $\psi(N)=N + \lambda$ ,
 $\lambda$ being in the field, we recover the commutation relations of the usual
  harmonic oscillator algebra.

\nn  Generalizing the pioneer work of Bargmann \cite{Bargmann} for the usual 
harmonic oscillator, we have studied in \cite {canada}, \cite{barg1} and \cite{barg2}
 the Bargmann representations of
  the DHOA, defined by  (\ref{a3}). 

\begin{definition} \label{espacebargmann}
A Bargmann Hilbert space is a space ${\cal S}$  of
   functions, holomorphic on a ring $D$ of the complex plane,
  the
 scalar product of which is written with a
 true integral on the form :
 \be
(g,f) = \int F(z \zz ) f(z)\overline{g(z)}dz d\zz
\label{prosca}
\ee

\nn A Bargmann representation of a Deformed Harmonic Oscillator
 Algebra
 is a representation  on a 
 Bargmann Hilbert space such as  the annihilation or the creation operator admits
 eigenvectors generating ${\cal S}$.
\end{definition}

\nn Let us stress that in this definition, we discard the occurrence of
 q-integration in the scalar product, contrarily to many authors, see in
  particular \cite{gray}\cite{chaichian}.

\nn In this paper, starting from the reciprocal point of view of that developed
 in  \cite {canada}, \cite{barg1} and \cite{barg2}, we construct
the DHOA defined by a given Bargmann
  representation. That is, a Bargmann Hilbert space ${\cal S}$ 
   being given, we look for DHOA that can be represented by operators acting
    on ${\cal S}$.

\nn In section \ref{representations}, we recall briefly the irreducible representations of the
 DHOA  on the
 basis of the eigenvectors of $N$ and we discuss the existence of coherent
  states, defined as the eigenstates of the operators $a$
   (or $a^\dagger$). 
In section \ref{Bargmann}, we summarize the study of Bargmann representations for the DHOA, 
done in \cite {canada}, \cite{barg1} and \cite{barg2} and set up the 
  relations between the weight function $F$ defining the scalar product and
   the function $\psi$ characterizing the DHOA. Section \ref{DHOAfromF}
    is devoted to the true subject of this paper.
	  We study the conditions of existence of a function
	    $\psi$ and consequently of a DHOA defined by (\ref{a3})
		 when we impose that the DHOA has a representation on
		  a given Bargmann Hilbert space. We prove that if
		   the weight function F fulfills sufficient and necessary conditions, the construction
	 can be performed.
	 Furthermore, we obtain a necessary condition on the 
	 function $\psi$
	 in order that the DHOA admits a Bargmann representation
	 when the domain of the coherent states is a true ring. 
 In section \ref{examples}, we treat several examples  illustrating
  our construction, showing in particular that the
  sufficient conditions obtained in Section 4, are not necessary.
 
\section{Representations}\label{representations}

\subsection{Eigenvectors of $N$}

\nn Let  $\mid 0>$ be  the  eigenvector of $N$ with eigenvalue
 $\mu$ that is 
  assumed to be zero to simplify the notations and is restored when necessary for
  discussions.
We built the 
 normalized vectors  $\mid n >$  

\be
\mid n> =\left\{
\begin{array}{ll}
\lambda _n a^{\dagger n}\mid 0>,& \quad n\in N^+\\
\lambda _n a^{-n}\mid 0>,& \quad n\in N^-
\end{array}
\right.
\ee

\nn with 

\be
\lambda_n^{-2} =\psi (n)!= \left\{\begin{array}{ll}
1, &\quad n=0\\
\prod_{i=1}^n \psi(i),& \quad n \in N^+ - \{0\}\\
 \prod_{i=0}^{n+1} \psi(i),& \quad n \in N^-
\end{array}
\right.
\ee

\nn $N^+$ and $N^-$ are the set of integers $\geq 0$ and $<0$.

 \nn The vectors $\mid n>$ are the eigenvectors of $N$ with eigenvalue $n$
  and span the Hilbert space $\cal{H}$. The condition that
   $<n\mid aa^\dagger \mid n> $ is strictly positive, 
    restricts  the spectrum of $N$.
    The elements of $SpN$ are the integers in any interval $[\nu,\nu^\prime[$
	 in which 
	 the function $\psi$ 
	  is finite and strictly positive. Eventually one of the edge or both can
	   be infinity. When one edge is finite, it is an integer when $\mu =0$,
	   and a zero of $\psi$.
 We thus get different types of representations \cite{nous1},
	   \cite{nous2},
	   \cite {quesne}, \cite{Kos} according as $\psi$ has no zero, one zero or more, 
the distance between two consecutive 
zeros being  integer.
		 The representations are defined by~:

\be
\left\{
\begin{array} {ll}
a^\dagger \mid n >  =& (\psi (n+1))^{1/2} \mid n+1 > \\ 
a\mid n >  =& (\psi (n))^{1/2} \mid n-1 > ,\quad n \in SpN\\ 
N \mid n >  =&  (n )\quad \mid n > 
 \end{array}
\right.
\ee

\nn When $\mu$ is kept different from zero, al the relations of this section
 remain valid provided we change $\psi(\rho)$ into $\psi(\mu +\rho).$
 
 \begin{proposition}\label{mu.0}
 By construction, the eigenvalue $\mu$ of the starting state $\mid~0~>$
 belongs always to the spectrum of $N$.
 \end{proposition} 
\nn The first step to build a Bargmann representation requires to study the
 coherent vectors that constitute the basis vectors of the representation.

\subsection{Coherent states} \label{coherent}

 \nn We call coherent states, the eigenvectors of the operator $a$ or $a^\dagger$. 

\nn The state $\mid z>= \sum_p c_p\mid p>$  is an eigenvector of the
annihilation operator $a$ if
 the coefficients $c_p$ verify the recursive relation 

\be
z c_p = \psi (p+1)^{1/2}c_{p+1}
\label{cp}
\ee

\nn  When the spectrum of $N$  is finite, we have proved that $a$ and $a^\dagger$
 have no eigenvectors, hence  no Bargmann representation exists.

\nn  When the spectrum of $N$ is no upper bounded, let us denote this spectrum
 $SpN = \lambda + N^+$,
with the convention  $\lambda + N^+ =Z$ when $\lambda = -\infty$.
It results from Proposition \ref{mu.0}, that $\lambda \in N^-+\{0\}$.
The  eigenvectors $\mid z>$ of $a$ take the form~:

\be
\mid z>=\sum_{n=-1}^{\lambda}z^{n}(\psi(n)!)^{1/2}\mid n>
+\sum_{n=0}^{\infty}
z^n (\psi(n)!)^{-1/2}\mid n>
\label{z}
\ee

\nn The domain $D$ of existence of the coherent states depends on the function
 $\psi$. Indeed, $\mid z>$ belongs to the Hilbert space spanned by the basis
  $\mid n>$ only if the series in the right hand side of (\ref{z}) are convergent
   in norm, see the detailed discussion in \cite{canada}, \cite{barg1}, \cite{barg2}. 

\nn An analogous reasoning holds when $Sp N$ is no lower
 bounded, for the eigenstates of 
 $a^\dagger$, it results from Proposition \ref{mu.0} that, in this case,
  $\lambda \in N^+-\{0\}$.
  We have proved that  the eigenstates of $a$ and $a^\dagger$ never
 coexist. Let us summarize the results :
   
   \begin{proposition} \label{psia}
   The eigenvectors of the annihilation operator of a Deformed Harmonic
    Oscillator Algebra exist if the function $\psi$ occurring in the relations
	 (\ref{a3}) belongs to two classes : 

\nn $\bullet$  $\psi $ is a  strictly positive function without singularity
 on the whole real axis,
 $SpN =Z$,
 and
 
\nn $lim_{p \rightarrow -\infty} \psi(p)^{1/2} <lim_{p \rightarrow \infty}
 \psi(p)^{1/2}$, then the domain of existence of the coherent states is :

\nn $D = \{z ;\quad lim_{p \rightarrow -\infty} \psi(p)^{1/2}
 <\mid z\mid <lim_{p \rightarrow \infty}
 \psi(p)^{1/2}\}$. 

\nn $\bullet$  It exists $\lambda \in N^- + \{0\}$ such
 that $\psi(\lambda)=0 $ and $\psi$ is strictly
 positive  without singularity when $x > \lambda$, $SpN = \lambda +N^+$ then : 

\nn $ D =\{z;\quad \mid z\mid < lim_{p \rightarrow \infty}
 \psi(p)^{1/2}\}$.
 \end{proposition}
 
 \nn and

  \begin{proposition} \label{psiacroix}
   The eigenvectors of the  creation operator of a Deformed Harmonic
    Oscillator Algebra exist if the function $\psi$ occurring in the relations
	 (\ref{a3}) belongs to two classes : 

\nn $\bullet$  $\psi $ is a  strictly positive function  without singularity on the whole real axis, 
 $SpN =Z$,
 and
 
\nn $lim_{p \rightarrow -\infty} \psi(p)^{1/2} > lim_{p \rightarrow \infty}
 \psi(p)^{1/2}$,

\nn then the domain of existence of the coherent states is

\nn $D = \{z ;\quad lim_{p \rightarrow \infty} \psi(p)^{1/2}
 <\mid z\mid <lim_{p \rightarrow - \infty}
 \psi(p)^{1/2}\}$. 

\nn $\bullet$  It exists $\lambda \in N^+-\{0\}$ such that $\psi(\lambda)=0 $
 and $\psi$ is strictly
 positive  without singularity when $x < \lambda$, then $SpN = \lambda +N^-$ and

\nn $ D =\{z;\quad \mid z\mid < lim_{p \rightarrow -\infty}
 \psi(p)^{1/2}\}$.
 \end{proposition}

\nn When the eigenvalue
 of $\mid 0 >$ is  $\mu \ne 0$, from Proposition \ref{mu.0},
  the point $\mu$, instead of the origin,  belongs to $SpN$ and the
 statements of the propositions must be changed accordingly. 

\begin{proposition}\label{creation.annihilation}
 The nature of the operators (to be  creation or annihilation operators)
 is inverted in the change :

\be
a= a^{\prime\dagger},\quad N= -N^\prime -1, \quad \psi (\rho) =\psi^\prime(-\rho).
\ee
\end{proposition}

\nn From this change, we deduce that Proposition \ref{psiacroix} can be obtained
 from Proposition \ref{psia}.
 In the following section, we set up the Bargmann representation when the eigenvectors
of the annihilation operator exist and using
 Proposition \ref{creation.annihilation}
we obtain the corresponding results when the eigenvectors of the creation operator exist.

\section{ Bargmann representation} \label{Bargmann}

\nn Let  $\mid f>$ be one state of $\cal{H}$ :

\be
\mid f>= \sum _{n\in Sp N}f_n \mid n> , \quad \sum _{n\in Sp N}
 \mid f_n \mid^2< \infty
\label{ff}
\ee

\nn with $SpN = \lambda + N^+$ and $\lambda \in N^- +\{0\}$.
 Following the construction \cite{Bargmann}, in the Bargmann representation any state
 $\mid f>$ of $\cal{H}$  is represented by the function of a complex variable $z$,
 $f(z) = < \zz \mid f >$ :

\be
f(z) = 
 \sum _{n \geq 0} z^n f_n \psi (n)!^{-1/2} +
  \sum _{n < 0}^{\lambda} z^n f_n (\psi (n)!)^{1/2},  \quad \sum _{n\geq \lambda}
 \mid f_n \mid^2< \infty
\label{f}
\ee

\nn where the variable $z$ belongs to the domain $D$ of definition
 of the eigenvectors of $a$. 

\nn $\bullet$ Let us summarize the results when the eigenvectors of $a$ exist :

\begin{proposition}\label{Z}
Let the function $\psi$ characterizing the Deformed Harmonic Oscillator Algebra 
 (\ref{a3}) belong 
 to the first class described in Proposition \ref{psia}.
 The space  $\cal{S}$  of the Bargmann  representation
is constituted with
   holomorphic functions in 
 $D = \{z;\quad lim_{p \rightarrow -\infty} \psi(p)^{1/2
 }<\mid z \mid < lim_{p \rightarrow \infty} \psi(p)^{1/2}\}$,

\nn the  Laurent expansions of which read (\ref{f}) with $\lambda = -\infty$. 
 
\nn In particular, when
  $ \lim_{p \rightarrow -\infty} \psi(p)^{1/2
 } = 0$, $D$ is a disk excluding  the origin  that is a
  essential singularity point.
\end{proposition}

\begin{proposition} \label {SpN=nu+N}
Let the function $\psi$ characterizing the Deformed Harmonic
 Oscillator Algebra  (\ref{a3}) belong 
 to the second class described in Proposition \ref{psia}. The space  $\cal{S}$ 
   of the Bargmann  representation
is the subspace of the space of
    functions holomorphic in 
 $D = \{z; \quad 0 <\mid z \mid < lim_{p \rightarrow \infty} \psi(p)^{1/2}\}$,
 
 \nn the  Laurent expansions of which read (\ref{f}). 

\nn The functions of $\cal{S}$ are of the form $z^{\lambda} g(z)$ where
 $\lambda$ is the lowest bound of the spectrum of $N$ and $g(z)$ is holomorphic in $D +\{0\}$.
The origin is  a pole of multiplicity lower 
or equal to $-\lambda$. In particular, when $\lambda =0$, the functions of
  $\cal{S}$
 are holomorphic at the origin.
\end{proposition}
	 
\nn In particular, to the basis vectors $\mid n> , n\in Sp N$ correspond
 the monomials~: 

\be
< \zz \mid n > = \left\{ \begin{array}{ll} 
z^n  (\psi (n)!)^{-1/2} , \quad &n \geq 0\\
z^n (\psi (n)!)^{1/2}, \quad &n <0
\end{array}
\right.
\label{n}
\ee

\nn $\bullet$ Let us summarize the results when the eigenvectors of $a^\dagger$ exist :

\begin{proposition}\label{Z+}
Let the function $\psi$ characterizing the Deformed Harmonic Oscillator Algebra 
 (\ref{a3}) belong 
 to the first class described in Proposition \ref{psiacroix}.
 The space  $\cal{S}$  of the Bargmann  representation
is constituted with
   holomorphic functions in 
 $D = \{z;\quad lim_{p \rightarrow +\infty} \psi(p)^{1/2
 }<\mid z \mid < lim_{p \rightarrow -\infty} \psi(p)^{1/2}\}$,

\nn The functions of $\cal{S}$ can be expanded in  Laurent series of the form :

\be
f(z) = 
 \sum _{n \geq 0}^{\lambda -1} z^{-n} f_n \psi (n)!^{1/2} +
  \sum _{n < 0} z^{-n} f_n (\psi (n)!)^{-1/2},  \quad \sum _{n \leq \lambda -1}
 \mid f_n \mid^2< \infty
\label{fcroix}
\ee 

\nn with $\lambda = +\infty$.
 
\nn In particular, when
  $ \lim_{p \rightarrow \infty} \psi(p)^{1/2
 } = 0$, $D$ is a disk excluding  the origin  that is a
  essential singularity point.
\end{proposition}

\begin{proposition} \label {SpN=nu+N+}
Let the function $\psi$ characterizing the Deformed Harmonic
 Oscillator Algebra  (\ref{a3}) belong 
 to the second class described in Proposition \ref{psiacroix}. The space  $\cal{S}$  of the Bargmann  representation
is the subspace of the space of
    functions holomorphic in 
 $D = \{z; \quad 0 <\mid z \mid < lim_{p \rightarrow \infty} \psi(p)^{1/2}\}$
  that read (\ref{fcroix}).

\nn The functions of $\cal{S}$ are of the form $z^{-\lambda+1} g(z)$ where
 $\lambda-1$ is the upper bound of the spectrum of $N$ and $g(z)$ is holomorphic in $D +\{0\}$.
The origin is  a pole of multiplicity lower 
or equal to $\lambda-1$. In particular, when $\lambda =1$, the functions of  $\cal{S}$
 are holomorphic at the origin.
\end{proposition}

\nn In particular, to the basis vectors $\mid n> , n\in Sp N$ correspond
 the monomials~: 

\be
< \zz \mid n > = \left\{ \begin{array}{ll} 
z^{-n}  (\psi (n)!)^{1/2} , \quad &n \geq 0\\
z^{-n} (\psi (n)!)^{-1/2}, \quad &n <0
\end{array}
\right.
\label{ncroix}
\ee

\nn  The function of $z$,
 $G(\zeta z)=<~\zz~\mid~\zeta~>$ corresponds to the coherent state $\mid \zeta>$.
  The holomorphic functions belonging to
  $\cal{S}$ have analytical properties strongly depending on $\psi$.
   Their growth on the edge of $ D$ is controlled by the growth of 
   $\mid G(z\zz)\mid ^{1/2}$.

\nn A  Bargmann representation exists if  we can obtain  a positive real function
 $F(x)$ such as 

\be
\int F(z \zz ) \mid \zz ><\zz \mid dz d\zz =1
\label{1}
\ee

\nn where the integration is extended to the whole complex plane and where
 $F(\mid z \mid ^2 )$ contains the characteristic function of the domain $D$
  of existence of the  coherent states. The existence of (\ref{1})
   ensures that the scalar product
  of the representation takes the form (\ref{prosca}).

 \nn We easily prove  that,

\begin{proposition}
 In  a Bargmann representation :

\nn $\bullet$ either $a^\dagger$
  is the multiplication by $z$,  $a$ the operator $z^{-1}\psi(zd/dz)$
   and $N$ the operator $zd/dz$

\nn $\bullet$ or $a$
  is the multiplication by $z$,  $a^\dagger$ the operator $z^{-1}\psi(-zd/dz)$
   and $N$ the operator $-zd/dz-1.$

\end{proposition}

\nn Let us introduce the Mellin transform $\hat{F}(\rho)$ of $F(x)$ :

\be
\hat {F}(\rho)=\int_0^{\infty} F(x) x^{\rho -1} dx 
\label{mel}
\ee

\nn From (\ref{n}) and (\ref{1}), we deduce that  $\hat{F}(\rho)$ exists on all the integers belonging to
the spectrum of $N$ and verify the
 following condition~:

\be
\hat{F}(n+1)= \left\{ \begin{array}{lll}
 \psi(n)!, & n\geq 0&\\
 (\psi(n)!)^{-1}, & n <0& , n \in Sp N
\end{array}
\right.
\label{moments}
\ee

\nn Let us remark that 

\be
\hat {F}(\rho) \leq \hat {F}(n)+\hat {F}(n+1), \quad n \leq Re \rho <n+1
\label{nn}
\ee

\nn because $F(x)$ is a positive function.
Therefore, the Mellin transform of $F$ exists for any $\rho$  such as
  $Re \rho \in [\lambda +1,+\infty[$.

\nn Formula (\ref{moments}) is equivalent to

\be
\hat{F}(n+1)=\psi (n) \hat{F}(n), \quad \mbox{with} \quad \hat{F}(1) =1
\label{M}
\ee

\nn which ensures that the operators $a^\dagger = z$, 
 $a = z^{-1}\psi(zd/dz)$ be adjoint on the basis $\mid n>$. 
 In \cite{barg1} and 
 \cite{barg2}, we have discussed  the interpolation of this equation, 
 the simplest one reads :

\begin{proposition}\label{annihilation}
 The function $\psi$ characterizing the DHOA 
defined in (\ref{a3}) and 
   the Mellin transform of the weight
 function $F$ are related by the following equation :

\be
\hat {F}(\rho+1)= \psi(\rho)\hat {F}(\rho), \quad \mbox{with} \quad \hat{F}(1) =1.
\label{mel1}
\ee

\nn when the coherent states are the eigenvectors of the annihilation operator $a$.
\end{proposition}

\nn We obtain an analogous proposition when the basis vectors of the
 Bargmann Hilbert space are the eigenvectors of $a^\dagger$. In this case, the
 annihilation operator is the multiplication by $z$.
  
 \begin{proposition} \label{creation}
 The function $\psi$ characterizing the DHOA 
defined in (\ref{a3}) and 
   the Mellin transform of the weight
 function $F$ are related by the following equation :

\be
\hat {F}(-\rho+1)= \psi(\rho+1)\hat {F}(-\rho), \quad \mbox{with} \quad \hat{F}(1) =1.
\label{mel2}
\ee

\nn when the coherent states are the eigenvectors of the creation
 operator $a^\dagger$.
\end{proposition}

\nn When the eigenvalue $\mu$ is different from zero, the only change
 in Equations (\ref{mel1}) and (\ref{mel2})
is $\psi(\rho) \rightarrow \psi(\mu +\rho)$.

\nn Proposition \ref{annihilation} corresponds to functions $\psi$ that are finite
 and positive
on a non upper bounded interval $]\lambda,+\infty[$ with $\lambda \in N^- +\{0\}$, 
as proved in 
Proposition \ref{psia}. The 
Mellin transform of $F$ given by (\ref{mel1}) is finite and
 positive on the same interval.
When $\lambda$ is finite, $\hat{F}(\lambda)$ is infinite. When  $\lambda$ is infinite, 
$\lim _{\rho \rightarrow -\infty }\frac{\hat{F}(\rho+1)}{\hat{F}(\rho)} <
\lim _{\rho \rightarrow +\infty }\frac{\hat{F}(\rho+1)}{\hat{F}(\rho)}. $

\nn Proposition \ref{creation} corresponds to functions $\psi$ that are finite and positive
on a non lower bounded interval $]-\infty,\lambda^\prime +1[$ with
 $\lambda^\prime \in N^+$ as proved in Proposition \ref{psiacroix}.
 The 
Mellin transform of $F$ given by (\ref{mel2}) is finite and positive on
 the  interval $]-\lambda^\prime, +\infty[$.
When $\lambda^\prime$ is finite, $\hat{F}(-\lambda^\prime)$ is infinite.
 When  $\lambda^\prime$ is infinite, 
$\lim _{\rho \rightarrow -\infty }\frac{\hat{F}(-\rho+1)}{\hat{F}(-\rho)} >
\lim _{\rho \rightarrow +\infty }\frac{\hat{F}(-\rho+1)}{\hat{F}(-\rho)}. $

\nn These two cases lead to the same conclusions :

\begin{proposition}\label{mellinF}
 The 
Mellin transform of the weight function of the Bargmann representation is finite and  
positive  on a non upper bounded interval of the real axis.

\nn $\bullet$ When the lower bound of the interval is finite, it  is an negative integer
 or zero
 and 
 $\hat{F}$ is infinite at this point.
 
 \nn $\bullet$ When the interval is not lower bounded, $\hat{F}$ fulfills the following
 condition :
 
 \be
 \lim _{\rho \rightarrow -\infty }\frac{\hat{F}(\rho+1)}{\hat{F}(\rho)} <
\lim _{\rho \rightarrow +\infty }\frac{\hat{F}(\rho+1)}{\hat{F}(\rho)} 
\ee 
 \end{proposition}
 
\nn Let remark that, while the spectrum of the energy operator $N$ is 
a non lower or non upper
bounded interval in $Z$, the definition domain of the Mellin transform of 
the weight function must be
a non upper bounded interval of the real axis containing the origin in order
that the construction be possible.    

\section{Deformed harmonic oscillator constructed from a weight function}
\label{DHOAfromF}

\nn The purpose of the two following sections is to construct a DHOA assumed to admit a given Bargmann representation. More precisely, we start
  with a Bargmann Hilbert space ${\cal S}$ of functions, holomorphic  on a ring
 $D$
   of the complex plane. The scalar
 product is written on the form (\ref{prosca}) in 
   terms of  a given  function $F$ positive on the interval
    $]\alpha, \beta[$.

\nn The construction can be summarize as follows : We look for a function $\psi$,
 we associate to this function the  DHAO defined by (\ref{a3}). Then, remains a consistency condition :
  we must prove that this algebra admits a
 Bargmann representation on ${\cal S}$, as defined in section \ref{Bargmann}. 

\nn Once the DHOA is constructed, the spectrum of $N$ is obtained and, eventually,
 the representation space must be restricted according to Propositions \ref{Z},
  \ref{SpN=nu+N}, \ref{Z+},
  \ref{SpN=nu+N+}.

\nn First, in order to define the function $\psi$ by applying (\ref{mel1}) 
or (\ref{mel2}), the Mellin transform 
 $\hat{F}$ of the given
 weight function $F$ must exist. Secondly, in order to define the basis vectors of 
the representation as the eigenvectors 
of the annihilation operator $a$ (or $a^\dagger$),
  the function  $\hat{F}$ must be define on a non upper bounded interval 
of the real axis,  due to the Proposition \ref{mellinF}.

 \nn Let us denote  $\hat{F}$ the Mellin transform of
  the weight function $F$  :

\be 
\hat{F}(\rho)=\int _{\alpha}^{\beta} F(x) x^{\rho -1} dx
\label{f1}
\ee

\nn From Proposition \ref{mellinF}, we obtain the following necessary conditions
 that must satisfy $\hat{F}$ in order that the construction of the DHOA
  be possible :
   
\begin{proposition}\label{Fmellin}
When a Bargmann Hilbert ${\cal S}$ space is given, a necessary condition 
in order that a DHOA admits
 a representation on ${\cal S}$ is that the Mellin transform of
  the weight function defining
 the scalar product (\ref{prosca}) exists on a non
  upper bounded interval of the real axis.
  
  \nn $\bullet$ When the lower bound of the 
  interval is finite, it must be a negative integer $\nu$. The functions of 
  the representation space read~:
  
  \be
  f(z) = \sum_{n\geq \nu-1} \frac{z^n f_n}{\hat{F}(n+1)^{\frac{1}{2}}},
    \quad \sum _{n\geq \nu-1}
 \mid f_n \mid^2< \infty
 \label{eq:f}
 \ee  
  
  \nn  $\bullet$ When the interval is not lower bounded, $\hat{F}$ must 
  fulfill the following
 condition~:
 
 \be
 \lim _{\rho \rightarrow -\infty }\frac{\hat{F}(\rho+1)}{\hat{F}(\rho)} <
\lim _{\rho \rightarrow +\infty }\frac{\hat{F}(\rho+1)}{\hat{F}(\rho)}. 
\ee 

\nn The functions of the representation space read (\ref{eq:f})
  with $\nu=-\infty$.   
\end{proposition}
 
\nn When the previous conditions are satisfied, 
 we are faced with two possibilities according as we choose that the basis vectors
 of the Bargmann representation are the eigenstates of the annihilation or of 
 the creation operator :
 
 \begin{proposition} \label{etats.propres.de.a}
When the Mellin transform of the weight function satisfies the necessary
 conditions of Proposition \ref{Fmellin}, the function :
 
\be
 \psi(\rho) = \frac{\hat {F}(\rho+1)} {\hat {F}(\rho)}
\label{eq.psi}
\ee

\nn put in the relations (\ref{a3}) defines a DHOA that can be represented
 on the Bargmann Hilbert space spanned with the eigenvectors of the
 annihilation operator, equipped with the scalar product (\ref{prosca}).
 
 \nn In this representation, the creation is the multiplication by $z$.
 \end{proposition}
 
  \begin{proposition} \label{etats.propres.de.acroix}
When the Mellin transform of the weight function satisfies the necessary
 conditions of Proposition \ref{Fmellin}, the function :
 
\be
 \psi(\rho) = \frac{\hat {F}(-\rho+2)} {\hat {F}(-\rho+1)}
\label{eq.psicroix}
\ee

\nn put in the relations (\ref{a3}) defines a DHOA that can be represented
 on the Bargmann Hilbert space spanned with the eigenvectors of the
 creation operator, equipped with the scalar product (\ref{prosca}).
 
 \nn In this representation, the annihilation is the multiplication by $z$.
 \end{proposition}

 \nn When the eigenvalue $\mu$ of the starting state is different from zero,
  the right hand side
of (\ref{eq.psi}) or (\ref{eq.psicroix}) is unchanged.
 The left hand side is replaced by $ \psi(\mu+\rho)$, the 
parameter $\mu$ of the representation is involved in the data of the weight
 function. 

\nn All the construction is done in terms of the Mellin transform 
of the weight function.
  It is worthwhile to note that the
 reproducing kernel $G(x)$ is always expressed in terms of $ \hat{F}$, by the
 same simple formula : 

\be
G(x) = \hat{F}(1)\sum _{n\geq \nu-1}\frac{x^n}{\hat{F}(n+1)}
\label{GF}
\ee 

\nn The main point is to prove that the DHOA so constructed admits a representation on the Bargmann Hilbert space ${\cal S}$, that is we have to establish the consistency of the construction.

\subsection{$D= D_{\alpha \beta}\equiv \{z; \quad 0 <\alpha< \mid z \mid ^2
 < \beta < +\infty \}$ }

\nn From (\ref{f1}), we obtain the following inequalities :

\be
\left\{
\begin{array}{ll}
\hat{F}(\rho) <& \alpha^{-\mu}\hat{F}(\rho +\mu)\\
\hat{F}(\rho) <& \beta^{\mu}\hat{F}(\rho -\mu)
\end{array}
\right. , \quad \forall \mu > 0
\label{ineq}
\ee

\nn that ensure that if $\hat{F}(\rho)$  diverges for one value of $ \rho$,
 it always diverges and that the Mellin transform of $F$ does not exist on
  the real axis. In the following we assume that $\hat{F}$ exists. Then the Mellin transform (\ref{f1}) exists for any $\rho \in R $ and
    is a strictly positive function. The functions $\psi$ defined by (\ref{eq.psi})
and (\ref{eq.psicroix})	 also are  strictly positive functions. Using these
 functions $\psi$ in
 (\ref{a3}), we obtain the corresponding DHOAs.

\nn Since  the function $\psi$ is strictly positive on $R$, 
the spectrum of
  $N$ is $Z$.

\nn The construction  is achieved if we prove that the coherent states of the
 so constructed DHOAs exist when $\alpha <\mid z\mid ^2<\beta$. We then have
  to determine the behaviors of the function $\psi$ at $\pm \infty$.

\nn Let us write (\ref{f1}) :

\be
\hat{F}(\rho)= \rho^{-1} \beta^\rho 
\int_{\left(\frac{\alpha}{\beta}\right )^\rho}^1 
F(\beta	 x^{\frac{1}{\rho}})dx
\label{f3}
\ee

\nn We easily see that $\hat{F}(\rho)\simeq F(\beta)\rho^{-1}\beta^\rho$ when
 $\rho$ goes to $+\infty$, if $F(\beta)$ is finite and different from zero. Therefore,
  the function $\psi(\rho)$ given in (\ref{eq.psi}) (resp.  (\ref{eq.psicroix}))
   goes to $\beta$ when $\rho$ goes to $+\infty$ (resp. $-\infty$) .

\nn When 
  $F^{(k)}(\beta) = 0, k= 0, \cdots, b$ and when $F^{(b+1)}(\beta)$ is finite and not equal to zero, 
   $\hat{F}(\rho)\simeq (-1)^{b+1} F^{(b+1)}(\beta)\rho^{-(b+2)}\beta^{b+1+\rho}$ when
 $\rho$ goes to $+\infty$. Then,  the function $\psi(\rho)$ given in (\ref{eq.psi}) (resp.  (\ref{eq.psicroix}))
   goes to $\beta$ when $\rho$ goes to $+\infty$ (resp. $-\infty$) .

\nn The same reasoning holds for the behavior at $-\infty$, indeed
 we now write
 (\ref{f1})~:

\be
\hat{F}(\rho)= \rho^{-1} \alpha^\rho \int^{\left (\frac{\beta}{\alpha}
\right )^\rho}_1 F(\alpha x^{\frac{1}{\rho}})dx
\label{f4}
\ee

\nn We get now that the limit at $-\infty$ :

\nn When $F(\alpha)$ is finite and different from zero, 
 $\hat{F}(\rho)\simeq F(\alpha)\rho^{-1}\alpha^\rho$.
 
\nn When  $F^{(k)}(\alpha) = 0, k= 0, \cdots, a$ and when $F^{(a+1)}(\alpha)$ is finite 
and different from zero, $\hat{F}(\rho)\simeq (-1)^{a+1} F^{(a+1)}(\alpha)
\rho^{-(a+2)}\alpha^{a+1+\rho}$.
 On these conditions, $\psi(\rho)$  given in (\ref{eq.psi})
  (resp.  (\ref{eq.psicroix}))
   goes to $\alpha$ when $\rho$ goes to $-\infty$ (resp. $+\infty$).
 
 \nn As $\alpha < \beta$, the necessary condition given in
  Proposition \ref{Fmellin} is fulfilled.
 Using the results of section \ref{Bargmann}, we prove that the coherent
 states, eigenstates of $ a$ or $a^\dagger$ according to the choice of the 
 characteristic function, are 
 defined when $\alpha <\mid z\mid ^2<\beta$. This completes the proof of the
  consistency of the reconstruction.

\nn We obtain the following sufficient conditions in order
 that the construction be consistent :

\begin{proposition} \label {prop1} 
 Let $F(x)$ be a function defined on the
 interval
 $]\alpha,\beta[$ such as : 
	
	\nn  $\bullet$  the Mellin transform of $F$ exists,
	 
  \nn  $\bullet$   $F(\alpha)$ is finite and different from zero or
 $F^{(a+1)}(\alpha)$ is finite and different from zero when 
     $F^{(k)}(\alpha) =0$, for $ k= 0,\cdots , a$,

\nn  $\bullet$  $F(\beta)$ is finite and different from zero or 
$F^{(b+1)}(\beta)$ is finite and different from zero when
  $F^{(l)}(\beta)=0$, for $l=0,\cdots, b$.
	     
 \nn One can construct two 
	   Deformed Harmonic
  Oscillator
  Algebras that admits a representation on a Hilbert space constituted by 
   functions, holomorphic in the ring
    $D_{\alpha \beta}=\{z;\quad \alpha^{\frac{1}{2}}<\mid z\mid < \beta^{\frac{1}{2}}\}$ and
	 equipped with the scalar product :
   \be
(g,f) = \int F(z \zz ) f(z)\overline{g(z)} \theta (z\zz - \alpha)
 \theta(\beta -z\zz) dz d\zz
\ee
 $\bullet$ One of this DHOA corresponds to the characteristic function (\ref{eq.psi}) and
 in the Bargmann representation, its creation operator is the multiplication
  by $z$.
  
  \nn  $\bullet$ The second  DHOA corresponds to the characteristic function
   (\ref{eq.psicroix}) and
 in the Bargmann representation, its annihilation operator is the multiplication
  by $z$. 
\end{proposition}

\nn If the edge conditions are not fulfilled the reconstruction
 may exist but we
 have to establish in each specific cases that the limits at infinities are
  the expected ones.

\nn Finally, let us remark that 

\be
\alpha \int^\beta_\alpha F(x) x^{\rho-1} dx
 \leq \int^\beta_\alpha F(x) x^{\rho} x dx 
 \leq \beta \int^\beta_\alpha F(x) x^{\rho-1}  dx
\label{f5}
\ee

\nn Obviously, these inequalities hold when $\alpha =0$ or $\beta$ infinite. We 
obtain :

\be
\alpha \leq \psi(\rho) \leq \beta
\label{f6}
\ee

\nn This constitutes a necessary condition that the 
function $\psi$  must fulfill  in order that the DHOA have a Bargmann representation :

\begin{proposition} \label{prop2}
Let the function $\psi$ characterizing the deformed harmonic oscillator algebra
(\ref{a3})   be strictly positive and such that its limits
 at infinities are finite and verify $ \psi(-\infty)
 < \psi(+\infty)$ (resp. $ \psi(-\infty)
 > \psi(+\infty)$).
    The coherent states, eigenstates of the annihilation (resp. creation)
	 operator,  exist
	 in the ring of the complex plane
\be
\begin{array}{lll}
    &D=\{z ; \quad \psi(-\infty)^{1/2}<\mid z\mid
 < \psi(+\infty)^{1/2}\}& \nonumber \\ 
 \mbox{(resp.}&
 D=\{z ; \quad \psi(+\infty)^{1/2}<\mid z\mid
 < \psi(-\infty)^{1/2}\}  \mbox{).}\nonumber
\end{array}
\nonumber
\ee
 The DHOA can be represented on a space of functions holomorphic in $D$ only if
  $\psi$ takes all its
    values  between the two limiting values :
\be 
\begin{array}{lll}
	 & \psi(-\infty)<\psi(\rho)
 < \psi(+\infty)& \nonumber \\
\mbox{(resp.}&
   \psi(+\infty)<\psi(\rho)
 < \psi(-\infty)&\mbox{).} \nonumber
\end{array}
\nonumber
\ee
	\end{proposition}

\subsection{$D= D_{0\beta} \equiv \{z; \quad  \mid z \mid ^2
 < \beta < +\infty \}$ } 

\nn The second inequality of (\ref{ineq}) still holds. Let us denote by
 $\nu$ the  value such as  the integration (\ref{f1}) is divergent when
  $ \rho \leq \nu$  and convergent when $ \rho > \nu$. When $\nu = -\infty$
   the Mellin transform exists for any $\rho$, while when $\nu = +\infty$ 
it never exists. As $\hat{F}(\nu)$ is infinite while $\hat{F}(\nu +1)$ is finite,
 $\psi(\nu)=0$ (resp. $\psi(1-\nu)=0$) defined by (\ref{eq.psi}) (resp. (\ref{eq.psicroix})
vanishes. The spectrum of $N$ is $\nu +N^+$ ( resp. $2-\nu + N^-$)
 when $\nu$ is finite and $Z$
  when $\nu = -\infty$. From Proposition
   \ref{Fmellin}, $\nu$ must
    belong to $N^-+\{0\}$. 
   
\nn The expressions of the coherent states  contain one or two infinite
 summations according to $\nu$ is finite or not.
 
\nn $\bullet$ In the case where
  $\nu = -\infty$, $Sp N = Z$, the limit at $-\infty$ must be done in each
   specific case and we treat an example in the subsection (8.1). This limit is positive or zero since $\hat{F}(\rho)$ is positive, it must be zero in order that the domain of existence of the coherent states be consistent with the definition of the given Bargmann Hilbert space.

\nn Let us summarize this result :

\begin{proposition} \label {prop3}
 Let $F(x)$ be a function defined on the
 interval
 $]0,\beta[$ such as :

	\nn $\bullet$  the Mellin transform of $F$ exists on the whole real axis,
	
 \nn  $\bullet$   $F(\beta)$ is finite and different from zero or  $F^{(b+1)}(\beta)$
is finite and different from zero when
  $F^{(l)}(\beta)=0$, for $l=0,\cdots, b$,
	    
 \nn   $\bullet$ and  
$ \lim_{\rho \rightarrow -\infty} \frac{\hat{F}(\rho+1)}{\hat{F}(\rho)} =0$.
  
 \nn One can construct two 
	   deformed harmonic
  oscillator
  algebras, corresponding to characteristic functions given in (\ref{eq.psi})
  and (\ref{eq.psicroix}), that admit a representation on a Hilbert space constituted by 
   functions, holomorphic in the ring
    $D = \{z ;\quad 0<\mid z\mid < \beta^{\frac{1}{2}}\}$ and
	 equipped with the scalar product :
  \be
(g,f) = \int F(z \zz ) f(z)\overline{g(z)}  \theta(\beta -z\zz) dz d\zz 
\label{scalbeta}
\ee 
\end{proposition}

\nn $\bullet$ In the case where $\nu$ is finite, it remains only one infinite
 summation in (\ref{z}), this summation is convergent if $\mid z \mid^2 <
 \lim_{\rho \rightarrow +\infty} \psi(\rho)$.  As $\beta$ is finite,
  this limit is obtained by the same reasoning as in the previous case.
  In this case, as the spectrum of $N$ is lower or upper bounded according as
we choose $\psi$ defined by (\ref{eq.psi}) or by (\ref{eq.psicroix}),
the Laurent expansions of  the functions belonging to
 the representation space contain terms in $z^n$ with $n \geq \nu$
  as results from Propositions
 \ref{SpN=nu+N} and \ref{SpN=nu+N+}. We then have :

\begin{proposition} \label {prop4} 
 Let $F(x)$ be a function defined on the
 interval
 $]0,\beta<+\infty[$ such as :

\nn  $\bullet$ the Mellin transform of $F$ only exists when $\rho > \nu$,
 $\nu$ finite and belonging to $N^- +\{0\}$.

\nn  $\bullet$   $F(\beta)$ is finite and different from zero or  $F^{(b+1)}(\beta)$
is finite and different from zero when
  $F^{(l)}(\beta)=0$, for $l=0,\cdots, b$.

\nn  One can construct two 
	   deformed harmonic
  oscillator
  algebras, corresponding to characteristic functions given in (\ref{eq.psi})
  and (\ref{eq.psicroix}), that admit a representation on a Hilbert space equipped with a scalar product 
(\ref{scalbeta}) and constituted with
   functions holomorphic in the ring
    $D = \{z; \quad 0<\mid z\mid < \beta^{\frac{1}{2}}\}$ restricted by the condition that the origin
 be a pole of 
multiplicity lower or equal to $\nu$.

\nn In particular when $\nu = 0$, the functions of the representation
space are holomorphic in a disk including the origin.

\end{proposition} 
 
\nn When $\beta$ is infinite, the consistency of the demonstration must be done 
in each case, one example is given in subsection 8.2.

\subsection{$D = D_{\alpha \infty} \equiv \{z; \quad 0< \alpha \mid z \mid ^2
 \} $}

\nn In this case, the integration (\ref{f1}) can diverge for $\nu< \rho$ and
 converge for $ \rho \leq \nu$. If $\nu$ is finite, $\psi(\nu)$ defined by equation
 (\ref{eq.psi}) and  $\psi(1-\nu)$ defined by equation
 (\ref{eq.psicroix})
is infinite. This corresponds to a case where $\psi$ has a singularity
    at a finite distance, not considered in this paper. If $\nu = +\infty$,
	 (\ref{f1}) always converges and the spectrum of $N$ is $Z$, the consistency
	  must be verified in each case.

\nn In this section, we have obtained consistency sufficient conditions to
 construct of a DHOA from its Bargmann representation when $D$ is a true ring
 or a true disk in the complex plane, an example will be given
 in the subsection 5.3. 

\section{Examples of construction}\label{examples}
 
\nn In this section, we give four examples of  construction of DHOA
 when its Bargmann representation is given, namely $F$ and 
 the domain $D$ of existence of the coherent states are given. In the first two
  examples $D$ is the whole complex plane, then the sufficient conditions of
consistency of the previous section do not apply. In the last two ones it is a
   ring,  one of them illustrates the results of the Propositions \ref{prop2}
   and\ref{prop3},
while in the last one  the sufficient conditions of the previous section are not fulfilled.
  
   \nn In the following,  the proofs are given for DHOA  resulting from 
   characteristic functions given in (\ref{eq.psi}) and for which the coherent
   states are the eigenvectors of the annihilation operator. Obviously, the
   same  can be developed when $\psi$ is given by (\ref{eq.psicroix}),
    leading to DHOA for which the coherent states are the eigenvectors
	 of the creation operator,  we just state the results.

\subsection{ $D= C-\{0\}$  and $  F(x)=\exp (-\sigma (\ln x)^{2n})$}

\nn  We assume that $ \sigma$ is a positive real number and that $ n$ is a positive integer.
As  the domain of existence of the coherent states is the whole
 complex plane, the sufficient condition of the previous section are not fulfilled.
 This example is an illustration of a case where the existence of the DHOA is established
 though the characteristic function $\psi$ is not obtained on an explicit form.

\nn The Mellin transform of $F(x)$ (\ref{mel}) reads :

\be 
\hat{F}(\rho) = \int_{-\infty}^{+\infty}e^{-\sigma t^{2n}+\rho t} dt
\label{mu}
\ee

\nn As $\sigma $ is positive and $n$ is a positive integer, we see
 that $\hat{F}(\rho)$ exists and is strictly positive for all $\rho$. The
  function $\psi$ given by (\ref{eq.psi}) is then a strictly positive function
   and the spectrum of $N$ is $Z$.

\nn The reconstruction of the deformed algebra
    will be achieved if we prove that the coherent states resulting of the
	 function $\psi$ such obtained are defined in the whole complex plane as
	  assumed. We thus have to study the behavior of $\psi$ at infinities.

\nn Let us assume that $\rho$ is positive. We write (\ref{mu}):

\be 
\hat{F}(\rho) = \int_{-\infty}^{+\infty}e^{-\sigma (t+v)^{2n}+\rho (t+v)} dt
\label{mu1}
\ee

\nn where we choose 

\be
v =\left(\frac{\rho}{\sigma(2n-1)}\right)^\frac{1}{2n-1}
\ee

\nn in order that the term under the exponential does not contain linear term in
 $t$ and (\ref{mu1}) reads :

\be
\hat{F}(\rho) = e^{-2(n-1)\sigma^\frac{1}{2n-1}
\left(\frac{\rho}{2n-1}\right)^\frac{2n}{2n-1}}
\int_{-\infty}^{+\infty}dt e^{-\sigma\left(\sum_{p\geq 2}C_p^{2n}t^p v^{2n-p}
\right)}
\label{mu2}
\ee

\nn After a change of variable $t=u v^{1-n}$, (\ref{mu2}) can be written :

\be
\hat{F}(\rho) = e^{-2(n-1)\sigma^\frac{1}{2n-1}\left(\frac{\rho}{2n-1}
\right)^\frac{2n}{2n-1}}\int_{-\infty}^{+\infty}du v^{1-n}
 e^{-\sigma n(2n-1) u^2 } e^{-\sigma\left(\sum_{p\geq 3}C_p^{2n}u^p v^{n(2-p)}
 \right)}
\label{mu3}
\ee

When $\rho \rightarrow +\infty$, the integral goes to zero like $v^{1-n}$ and
using (\ref{ff}), we have~:

\be
\psi(\rho)\simeq e^{\frac{4n(n-1)}{(2n-1)^2}\left(\frac{\rho}{\sigma(2n-1)}
\right)^\frac{1}{2n-1}}
\ee

\nn We therefore get :
\be
\lim_{\rho \rightarrow +\infty} \psi(\rho) = +\infty
\ee
 
\nn From (\ref{mu}), we deduce that
\be
\hat{F}(-\rho)=-\hat{F}(\rho)
\ee

\nn so that 

\be
\psi(-\rho)=\frac{\hat{F}(-\rho+1)}{\hat{F}(-\rho)}=\frac{\hat{F}(\rho-1)}
{\hat{F}(\rho)}=\frac{1}{\psi(\rho-1)}
\ee

\nn Thus when $\rho \rightarrow -\infty$, we get that

\be
\lim_{\rho \rightarrow -\infty} \psi(\rho)=0. 
\ee

\nn The domain of
 existence of the coherent states is the complex plane and the consistency of
  the reconstruction is established.

\nn As $SpN = Z$,  the representation space  is the space of  the functions holomorphic 
in the complex plane without the origin that is a essential singularity point. Similarly,
to Equation (\ref{eq.psicroix}) corresponds another DHOA.

\begin{proposition}
One can construct two Deformed Harmonic Oscillator Algebras that can be represented on the Bargmann Hilbert
space of functions holomorphic in the complex plane without the origin equipped with the scalar product :

$$(g,f) = \int \exp (-\sigma (\log z\zz )^{2n}) \overline{g(z)} f(z) dz d\zz $$

\nn $\bullet$ The characteristic function $\psi$ involved in (\ref{a3}) is :

$$\psi (\rho) =\frac{\int^{+\infty}_{-\infty} \exp (-\sigma t^{2n}+(\rho +1) t ) dt}
{\int^{+\infty}_{-\infty} \exp (-\sigma t^{2n}+\rho t ) dt}$$

 \nn and the spectrum of $N$ is $Z$ and the coherent states are the eigenvectors of $a$.

\nn  $\bullet$ The characteristic function $\psi$ involved in (\ref{a3}) is :

$$\psi (\rho) =\frac{\int^{+\infty}_{-\infty} \exp (-\sigma t^{2n}+(2-\rho ) t ) dt}
{\int^{+\infty}_{-\infty} \exp (-\sigma t^{2n}+(1-\rho) t ) dt}$$

 \nn and the spectrum of $N$ is $Z$ and the coherent states are the eigenvectors
 of $a^\dagger$.
\end{proposition}

\nn In the next subsection, we give another example where the results of the
 previous section do not apply and in which the explicit calculation of the characteristic
 function $\psi$ can be
 performed. The main interest of the next example is to be a generalization of 
the usual harmonic oscillator algebra.

\subsection{$D= C - \{0\}$ and $ F(x)=\exp (-  x^{\frac{k}{m}})$ }

\nn We assume that $\frac{k}{m}$ is put on an irreducible form and that it is
 positive. When $\frac{k}{m}=1$, $F$ is the weight function of the Bargmann
  representation of the usual harmonic oscillator \cite{Bargmann}.

\nn The Mellin transform of $F$ (\ref{mel}) reads :

\be
\hat{F}(\rho)=\int_0^{+\infty}e^{-x^{\frac{k}{m}}}x^{\rho-1} dx
\ee

\nn After a change of variable $u= x^{\frac{k}{m}}$, it reads 

\be
\hat{F}(\rho)=\frac{m}{k}\int_0^{+\infty}e^{-u}u^{\rho\frac{m}{k}-1} du = 
\frac{m}{k} \Gamma (\rho \frac{m}{k})
\ee

\nn and the function $\psi$ characterizing the DHOA and resulting from
 (\ref{eq.psi}) is  :

\be
\psi(\rho)= \frac{\Gamma(\frac{m}{k}(\rho +1))}{\Gamma(\frac{m}{k}\rho)}
\label{k/m}
\ee

\nn From this explicit expression of $\psi$, we deduce that this function is
 strictly positive on the positive  axis and vanishes at the origin. The
  spectrum of $N$ is $N^+$.

\nn Using the asymptotic behavior of $\Gamma(z)$ for large values of
 $\mid z \mid$, we get :

\be
\psi(\rho)\simeq \left( \frac{m}{k}\rho \right)^{\frac{m}{k}}
\ee

\nn Thus $\lim_{\rho\rightarrow +\infty} \psi(\rho)$ is infinite and the
 coherent states, as assumed, are defined in the whole $C$. As $Sp N = N^+$, the functions
of the representation space are holomorphic in $C$, including the origin. A similar
 construction can be performed with the characteristic function (\ref{eq.psicroix}).

 \begin{proposition}
One can construct two Deformed Harmonic Oscillator Algebra that can be represented
 on the Bargmann Hilbert
space of functions holomorphic in the whole complex plane  equipped with the scalar product :

$$(g,f) = \int \exp (- ( z\zz )^{\frac{k}{m}}) \overline{g(z)} f(z) dz d\zz ,
 \quad \frac{k}{m} >0  $$

\nn $\bullet$ The characteristic function $\psi$ is :

$$\psi (\rho) =\frac{\Gamma( \frac{m}{k}(\rho +1))}
{\Gamma( \frac{m}{k}(\rho))}$$

\nn The spectrum of $N$ is $N^+$ and the coherent states are the eigenvectors of $a$.

\nn $\bullet$ The characteristic function $\psi$ is :

$$\psi (\rho) =\frac{\Gamma( \frac{m}{k}(2-\rho))}
{\Gamma( \frac{m}{k}(1-\rho))}$$

\nn The spectrum of $N$ is $N^- +\{0\}$ and the coherent states are the eigenvectors of
 $a^\dagger$.

 \end{proposition}

\nn When $\frac{k}{m}=1$, the function $\psi$ resulting of Equation (\ref{k/m}) is
the characteristic function of the usual harmonic oscillator in a fixed
 representation, namely $a^\dagger a = N$.
 
\nn We end this subsection by comparing the Bargmann representation considered
 in this subsection with the Bargmann representation of the usual harmonic
  oscillator.

\nn In this subsection, the scalar product (\ref{prosca}) is defined in the space
 $\cal{S}$ of  holomorphic functions of one complex variable and reads :

\be
(g,f)=\int d\zeta d\overline{\zeta} e^{-\zeta\overline{\zeta}^\frac{k}{m}}
 f(\zeta)\overline{g(\zeta)} , f , g \in \cal{S}
\label{u3}
\ee

\nn Denoting $\zeta=\chi e^{i\tau}$, it can be written :

\be
(g,f)= m^{-1}\int_0^{2\pi m} \frac{d\tau}{2\pi}\int _0^{+\infty}
 d\chi^2 e^{-\chi^\frac{2 k}{m}} f(\chi e^{i\tau})\overline{g(\chi e^{i\tau})}
\label{u4}
\ee

\nn The scalar product for the usual Bargmann representation reads :

\be
(g_B,f_B)=\int dzd\zz e^{-z\zz} f_B(z)\overline{g_B(z)} ,  f_B , g_B\in
 \cal{S}_B
\label{u}
\ee

\nn It takes the form :

\be
(g_B,f_B)= k^{-1}\int_0^{2\pi k} \frac{d\theta}{2\pi}\int _0^{+\infty}
 d\rho^2 e^{-\rho^2} f_B(\rho e^{i\theta})\overline{g_B(\rho e^{i\theta})}
\label{u1}
\ee

\nn Let us change $z=\zeta^{ \frac{k}{m}}$, we see that $0\leq \tau < 2\pi m$
 and that (\ref{u1}) reads :

\be
(g_B,f_B)= \frac{k}{m^2} \int_0^{2\pi m} \frac{d\tau}{2\pi}\int _0^{+\infty}
 d\chi^2 e^{-\chi^{ \frac{2 k}{m}}} \zeta^{ \frac{k}{m}-1}f_B(\zeta ^\frac{k}{m})
 \overline{\zeta^{ \frac{k}{m}-1}g_B(\zeta^\frac{k}{m})}
\label{u2}
\ee

\nn The scalar products written in (\ref{u2}) and in (\ref{u4}) are the same but
 they are not defined on the same space of functions :

\nn Indeed let us write $f_B(z) =\sum _{l\leq 0}f_l z^l$, the functions 
$f(\zeta) \equiv \zeta^{ \frac{k}{m}-1}f_B(\zeta ^\frac{k}{m})$ belong to
 $\cal{S}$ iff $f_l =0$ when $l \not= nm -1$, $n$ being a strictly positive
  integer. The functions $f$ such obtained belong to $\cal{S}$ but do not cover
   the whole space for they read :

\be
f(z) = \sum_{n=1}^{+\infty}f_{nm-1} z^{kn-1}
\label{u5}
\ee

 \begin{proposition}
Let us consider the two Bargmann Hilbert spaces on which are represented 
 the usual harmonic oscillator algebra and
 the DHOA considered in this subsection.  When, by 
a change of variables, their
 scalar products are written on the
  same form (\ref{u4}), the functions
	 belonging to the intersection of these two spaces are of the form (\ref{u5}).
\end{proposition} 

\subsection{$D = D_{\alpha \beta} 
  $ and
   $ F(x)= x^\sigma $}
 
\nn We start with a Bargmann representation such as the coherent states are
 defined on a ring of the complex plane $0 \leq \alpha \leq \rho \leq \beta
  < +\infty$.  This subsection illustrates the previous section with an example
 where we obtain an explicit expression for the characteristic function $\psi$. 

\nn The Mellin transform of the weight function reads :

\be 
\hat{F}(\rho)=\int _\alpha ^\beta x^{\sigma+\rho -1} dx
\label{x}
\ee 

\nn First, we see that this integration is finite for any $\rho$ and any $\sigma$
 when $\alpha \not= 0$ and for $\rho > -\sigma$ when $\alpha=0$.

\nn The resulting function $\psi$, defined by (\ref{eq.psi}), takes the form :

\be
\psi (\rho ) = \frac{\sigma +\rho}{\sigma +\rho +1 } \frac{\beta^{\sigma +\rho +1 }
 -\alpha^{\sigma +\rho +1 }}{\beta^{\sigma +\rho } -\alpha^{\sigma +\rho  }}
\label{x1}
\ee

\nn We now must look for the domain of existence of the coherent states in order
 to verify the consistency of this construction.

\nn $\bullet$ When $\alpha \ne 0$, the function $\psi$ is always positive and the spectrum
of $N$ is $Z$. It is easy to find that  the
  function $\psi(\rho)$ goes to $\alpha$ or $\beta$ when 
  $\rho \rightarrow -\infty$ or $+\infty$. This implies that the coherent states
   are defined for $\alpha \leq \rho^2 \leq \beta$, as expected. As $SpN = Z$, no
 restrictions appear on the Laurent expansions of the holomorphic functions of the
representation space. The same construction can be done starting with the characteristic
function (\ref{eq.psicroix})

\begin{proposition}
One can construct two Deformed Harmonic Oscillator Algebras that can be represented on
the space of functions holomorphic in $D_{\alpha \beta} = \{z;\quad 0 < \alpha < \mid z \mid ^2 < \beta\}$
equipped with 
the scalar product :

$$(g,f) = \int_{ 0 < \alpha <\mid z\mid^2<\beta}
 (z \zz)^\sigma \overline{g(z)} f(z) dz d\zz,
  \quad \forall \sigma. $$

\nn The characteristic functions are $\psi (\rho )$ or $\psi (1-\rho )$ expressed in 
(\ref{x1}), $SpN$ is $Z$ in both cases and the coherent states
are the eigenvectors of the annihilation or of the creation operator.
 
\end{proposition}

\nn $\bullet$ When $\alpha = 0$, the characteristic function resulting from (\ref{eq.psi}) is :

\be
\psi (\rho) = \frac{\sigma +\rho}{\sigma +\rho +1}\beta
\label{alpha=0}
\ee
From Proposition \ref{mellinF}, we deduce that the construction is only possible when  
$\sigma \in N^+$. Then $SpN = -\sigma +N^+$. The coherent states are defined in 
 $D = \{z;\quad \mid z \mid ^2 < \beta\}$ as expected and the origin is a pole of multiplicity
lower than $\sigma$ for the functions of the representation
space. Let us summarize this result and that obtained starting with Equation
 (\ref{eq.psicroix}) :

\begin{proposition}\label{F=xsigma}
One can construct two Deformed Harmonic Oscillator Algebras  that can be represented on
the space of functions of the form $z^{-\sigma} f_0(z)$ where $f_0(z)$ is 
 holomorphic in the whole disk  $D_{0\beta} = \{z;\quad \mid z \mid ^2 < \beta\}$ equipped with  
the scalar product :

\be
(g,f)=\int_{\mid z\mid^2<\beta} (z\zz)^{\sigma} f(z)
 \overline{g(z)}dz d\zz
\ee
provided that $\sigma$ be zero or a positive integer.
Their characteristic functions are $\psi (\rho)$ or $\psi (1-\rho)$ written in Equation
(\ref{alpha=0}).
 The spectrum of $N$ is $-\sigma + N^+$ or $\sigma + N^- +\{0\}$ and the coherent states
are the eigenvectors of the annihilation or of the creation operator. 
\end{proposition}

\subsection{$D = D_{0 \beta} 
  $ and
   $ F(x)= x^\sigma (\beta - x)^\eta $}

\nn This example does not fulfill the general conditions of Section \ref{DHOAfromF} for
the derivatives of the weight function on the edge $\beta$ of $D$ is zero or infinite when
$\eta$ is not a positive integer. The Mellin transform of $F$ exists when $\rho +\sigma >0$
and $\eta +1 >0$ and can be calculated in terms of
the $B$ function \cite{grad} :

\be
\begin{array}{ll}
\hat{F}(\rho)&=\int _0 ^\beta (\beta-x)^\eta x^{\sigma+\rho -1} dx\\
&=\beta^{\eta+\sigma+\rho}B(\rho+\sigma,\eta+1)\\
&=\beta^{\eta+\sigma+\rho}\frac{\Gamma(\rho+\sigma) \Gamma(\eta +1)}
{\Gamma(\rho+\sigma+\eta +1)}\\
\end{array}
\ee 
The expression of $\hat{F}(\rho)$ put in Equation (\ref{eq.psi}) leads to :

\be
\psi(\rho)= \frac{\rho+\sigma}{\rho+\sigma+\eta+1}\beta
\label{zerobeta}
\ee
 When $\eta =0$, one recovers the function (\ref{alpha=0}). The reasoning
  and the results
are similar to those of Proposition \ref{F=xsigma} :

\begin{proposition}
One can construct two Deformed Harmonic Oscillator Algebras  that can be represented on
the space of functions of the form $z^{-\sigma} f_0(z)$ where $f_0(z)$ is 
 holomorphic in the whole disk  $D = \{z;\quad \mid z \mid ^2 < \beta\}$ equipped with  
the scalar product  :

\be
(g,f)=\int_{\mid z\mid^2<\beta}(\beta - z\zz)^\eta (z\zz)^{\sigma}f(z)
 \overline{g(z)}dz d\zz
\ee
provided that $\sigma$ be zero or a positive integer.
Their characteristic functions  are $\psi (\rho)$ or $\psi (1-\rho)$ written in Equation 
(\ref{zerobeta}).
 The spectrum of $N$ is $-\sigma + N^+$ or $\sigma + N^- +\{0\}$ and the coherent states
are the eigenvectors of the annihilation or of the creation operator. 
\end{proposition}

\subsection{$D = D_{\alpha \beta} 
  $ and
   $ F(x)= \exp(x-\beta)^{-1} $}
   In this example, the general conditions of the section \ref{DHOAfromF} are
    not fulfilled too, for all the derivatives of $F(x)$ vanish on the edge
	 $\beta$ of $D$. To simplify,
 we give the proofs for
 $\alpha = 0$.  

\nn The Mellin transform of $F(x)$ are defined for $\rho >0$.
	 Integrating by parts, we get :
	 \be
	 \hat{F}(\rho) = \frac{1}{\rho}\int _{0}^{\beta} exp\left(\frac{1}{x-\beta}\right)
	 \frac{x^\rho}{(\beta-x)^2} dx
	 \ee
	 Expanding $(\beta-x)^{-2}$, we obtain the relation :
	 
	 \be
	 \rho = \beta^{-2}\sum_{n\geq 0} \frac{n}{\beta^n}\frac{\hat{F}(\rho +n+1)}  
 {\hat{F}(\rho )}
 \label{ro}
 \ee
 Now, when $\rho \rightarrow \infty$, $\frac{\hat{F}(\rho+1)}{\hat{F}(\rho)}$
necessarily
  goes to the limit $\beta_0 \leq \beta$ according Proposition \ref{prop2}. This
  implies :
\be
\lim_{\rho \rightarrow \infty}\frac{\hat{F}(\rho+n+1)}{\hat{F}(\rho)} =\beta_0^n 
 \ee
This result is used to calculate the right hand side of Equation (\ref{ro})
 when $\rho$ goes to $ \infty$. As the limit of the left hand side is infinite, we
  obtain
   a contradiction unless $\beta_0 \ne \beta$.
  As $\hat{F}(0) $ is infinite, $\psi(0)=0$ and the coherent states are defined 
  in $D_{0 \beta}$. The consistency conditions are satisfied. Let us 
state the result for arbitrary $\alpha$ :

\begin{proposition}
One can construct two Deformed Harmonic Oscillator Algebras  that can be represented on
the space of functions 
 holomorphic in   $D_{\alpha \beta}$  equipped with  
the scalar product  :

\be
(g,f)=\int_{ \alpha <\mid z\mid^2<\beta} \exp(z\zz-\beta) f(z)
 \overline{g(z)}dz d\zz
\ee
 The spectrum of $N$ is $ N^+$ or $N^- +\{0\}$ when $\alpha=0$ and $Z$ when
 $\alpha \ne 0$ and the coherent states
are the eigenvectors

\nn $\bullet$ of the annihilation  when the 
characteristic function is
 
\be
\psi (\rho) = \frac{\int _{\alpha}^{\beta} \exp(x-\beta)^{-1}
	 x^{\rho} dx} {\int _{\alpha}^{\beta} \exp(x-\beta)^{-1}
	 x^{\rho-1} dx}
\ee
\nn $\bullet$ or of the creation operator when the characteristic function is $\psi(1-\rho)$.
 
\end{proposition}
   
   \nn From these  examples, we conjecture that the propositions of
   Section \ref{DHOAfromF} can be largely extended.
    
\section{Conclusion}

\nn Giving a ring $D$ in the complex plane and a
  positive function $F$  that characterize a functional Bargmann Hilbert space $\cal{S}$,
we have discussed  the conditions under which exists a DHOA that 
admits a Bargmann representation in  $\cal{S}$, as defined
in \ref{espacebargmann}. 
We have obtained conditions on the weight function in order that solutions  exist :

\nn - in Proposition \ref{Fmellin}, we give necessary conditions,

\nn - in Propositions \ref{prop1}, \ref{prop2}, \ref{prop3}, \ref{prop4}, we give
 sufficient conditions, when D is not
 the whole complex plane.

\nn When one DHOA solution exists, another DHOA exists, if the annihilation operator
of the first one possesses eigenvectors generating the representation space, the same 
holds for the creation operator of the second one.

\nn Finally, we have developed some examples that does not fulfill the sufficient
conditions, in particular, when $D$ is 
the whole complex plane. We have obtained  deformations of usual Harmonic Oscillator Algebra
through deformations of its  Bargmann representation.

\end{document}